\documentclass[11pt,draft]{article}

\usepackage{amsfonts,latexsym,amsthm}

\newenvironment{cases}{\left\{\begin{array}{cl}}{\end{array}\right.}
\def\txt#1{\quad\mbox{#1}\quad}

\newtheorem*{proposition}{Proposition}
\newtheorem*{theorem}{Theorem}
\theoremstyle{definition}
\newtheorem*{definition}{Definition}
\def\pf{\ifvmode\else\newline\fi\noindent\textsc{Proof:\ }}
\def\qed{\mbox{ $\Box$}}

\def\opC{{\mbox{\boldmath$\mathfrak{C}$}}}

\def\colon{\;:\;}
\def\tfrac#1#2{{\textstyle\frac{#1}{#2}}}

\begin{document}
\def\footnotemark{}
\author{J\"org Wenzel}
\title{A supplement to my paper on\\real and complex operator
  ideals\thanks{Research supported by German Academic Exchange Service
  (DAAD)}}
\maketitle
\begin{abstract}
  We show that the main problem left open in \cite{wen95} can be
  solved using the Banach spaces $Z_\alpha$ recently constructed by
  Kalton \cite{kal95}. This gives an example of a complex operator
  ideal that has no real analogue. It thus
  shows the richer structure of complex operator ideals compared with
  the real ones.
\end{abstract}

\section{Introduction}
\label{sec:intro}

For a Banach space $X$ let $\overline X$ denote its complex conjugate,
i.e.~the Banach space $X$ equipped with the scalar multiplication
$\alpha \odot x:= \overline \alpha x$.

For an operator $T:X\to Y$ let $\overline T$ denote its complex
conjugate, i.e.~the operator $T$ acting between $\overline X$ and
$\overline Y$. This is again a linear operator.

For a complex operator ideal $\opC$ the complex conjugate ideal
$\overline \opC$ is defined as the ideal consisting of all $T$ such
that $\overline T\in \opC(\overline X,\overline Y)$.


The ideal $\opC$ is called self conjugate if $\opC=\overline \opC$.

The problem of the existence of operator ideals that are not self
conjugate was left open in \cite{wen95}. However, it was reduced there
to the problem of finding a complex Banach space that is not
isomorphic to its complex conjugate, but that is isomorphic to its
square.

In \cite{kal95}, Kalton gives examples of Banach spaces $Z_\alpha$
which are not isomorphic to their complex conjugates. On the other
hand, the elementary nature of these spaces makes it possible to
easily verify that they are isomorphic to their squares, thus
yielding the required Banach spaces.

\section{The spaces $Z_\alpha$}
\label{sec:space}

In \cite{kal95} Kalton constructed elementary examples of Banach
spaces $Z_\alpha$, not isomorphic to their complex conjugates. We
need only one additional property of these spaces which is not
explicitly mentioned in \cite{kal95}, namely that $Z_\alpha$ is
Cartesian. This is however easy, as is shown below.

For the convenience of the reader, let us repeat here the definition
of $Z_\alpha$.

Let $f_\alpha(t):= t^{1+i\alpha}$ for $-\infty<\alpha<\infty$ and
$0\leq t<\infty$. Given a sequence $x=(\xi_k)\in l_2$ define the
sequence $\Omega_\alpha(x)$ by
\[ (\Omega_\alpha(x))_k :=
\begin{cases}
  \xi_k f_\alpha\left( \log \tfrac{\|x\|_2}{|\xi_k|} \right) & \txt{if
    $\xi_k\not=0$} \\
  0 & \txt{otherwise.}
\end{cases}
\] Let $Z_\alpha$ be the space of all pairs of complex--valued
sequences $(x,y)$ such that
\[ \|(x,y)\|_\alpha := \|x\|_2 + \|y-\Omega_\alpha(x)\|_2 < \infty.
\] It turns out
that $Z_\alpha$ is a Banach space under a norm equivalent to the
quasi norm $\|\cdot\|_\alpha$.

Moreover $\overline{Z_\alpha}=Z_{-\alpha}$ and $Z_\alpha$ is
isomorphic to $Z_\beta$ if and only if $\alpha=\beta$.

\begin{definition}
  A Banach space $X$ is called \emph{Cartesian} if it is
  isomorphic to its Cartesian square $X\oplus X$.
\end{definition}

\begin{proposition} \label{cor}
  The spaces $Z_\alpha$ are Cartesian.
\end{proposition}
\pf
For a complex--valued sequence $x=(\xi_k)$ define
\[ U^{odd}x  := (\xi_{2k-1}) \quad \mbox{and} \quad
   U^{even}x := (\xi_{2k}).
\] Then the map $U:Z_\alpha \to Z_\alpha \oplus Z_\alpha$ defined by
\[ U(x,y) := (U^{odd}x,U^{odd}y) \oplus (U^{even}x,U^{even}y)
\] defines an isomorphism.

The bijectivity of $U$ is trivial, so we only show its continuity.

It is mentioned in \cite{kal95} that there is a constant $C$, such
that for $s\in l_\infty$ we have the estimate
\[ \|(sx,sy)\|_\alpha \leq C \cdot \|s\|_\infty \cdot
\|(x,y)\|_\alpha.
\] Hence
\[ \|(U^{odd}x,U^{odd}y)\|_\alpha \leq C
\|(x,y)\|_\alpha
\quad \mbox{and} \quad
\|(U^{even}x,U^{even}y)\|_\alpha \leq C
\|(x,y)\|_\alpha.
\] This proves the assertion. \qed

\section{Main theorem}
\label{sec:main}

\begin{theorem} \label{th:main}
  For $\alpha\not=0$ the operator ideal
  \[ \opC_\alpha := \{T\colon \mbox{\rm $T$ admits a factorization
    over the space $Z_\alpha$} \}
  \] is not self conjugate.
\end{theorem}
\pf Since $Z_\alpha$ is Cartesian, $\opC_\alpha$ is indeed an operator
ideal. Now $Id_{Z_\alpha}\in\opC_\alpha$. Assume that
$\overline{Id_{Z_\alpha}}=Id_{\overline{Z_\alpha}}\in\opC_\alpha$.
This implies that $Z_\alpha$ and $\overline{Z_\alpha}$ are
complemented in each other.  Since $Z_\alpha$ and
$\overline{Z_\alpha}$ are Cartesian, we can apply Pe\l{}czy\'nski's
decomposition method to obtain that $Z_\alpha$ and
$\overline{Z_\alpha}$ are actually isomorphic which contradicts the
properties of these spaces.  See \cite{wen95} for details.  \qed

\vspace{1ex}
\noindent
{\sc Mathematisches Institut, FSU Jena, 07740 Jena, Germany}\\
{\it E--mail:} {\tt wenzel@minet.uni-jena.de}


\begin{thebibliography}{1}

\bibitem{kal95}
{\sc Kalton, N.~J.}
\newblock An elementary example of a {Banach} space not isomorphic to its
  complex conjugate.
\newblock {\em Canad. Math. Bull.} {\bf 38} no.~2 (1995), 218--222.

\bibitem{wen95}
{\sc Wenzel, J.}
\newblock Real and complex operator ideals.
\newblock {\em Quaestiones Mathematicae} {\bf 18} no.~1--3 (1995), 271--285.

\end{thebibliography}
\end{document}